# SUBJECTIVE QUESTIONS AND ANSWERS FOR A MATHEMATICS INSTRUCTOR OF HIGHER EDUCATION

Florentin Smarandache
University of New Mexico
200 College Road
Gallup, NM 87301, USA
E-mail: smarand@unm.edu

**Abstract:** This article of mathematical education reflects author's experience with job applications and teaching methods and procedures to employ in the American Higher Education. It is organized as a standard questionnaire.

1) **What are the instructor's general responsibilities?**

   - participation in committee work and planning
   - research and innovation
   - on job training
   - participation to meetings
   - order necessary textbooks, audio-visual, and other instructional equipment for assigned courses
   - submit requests for supplies, equipment, and budgetary items in good order and on time
   - to keep abreast of developments in subject field content and methods of instruction
   - to assess and evaluate individual student progress, to maintain student records, and refer students to other appropriate college staff as necessary
   - participate on college-wide registration and advising process
   - effective and full use of the designated class meeting time
   - adequate preparation for course instruction, course and curriculum planning
   - teaching, advising students
   - be able to make decisions
   - knowledge and use of material
   - positive relationships
   - knowledge of content
   - plan and implement proposed plans (or change them if they don't work) – short and long-term planning
   - be a facilitator, motivator, model, appraiser and assessor of learning, counselor, classroom manager (i.e. manage student's behavior, the environment, the curriculum)
   - knowledge of teenage growth and development
   - continuously develop instructional skills.

The most important personal and academic characteristics of a teacher of higher education are: to be very good professionally in his/her field, to improve permanently his/her skills, to be dedicated to his/her work, to understand the students' psychology, to be a good educator, to deliver attractive and interesting lessons, to make students learn to think (to solve not only mathematical problems, but also the life ones), to try to approach mathematics with what students are good at (telling them, for example, that mathematics are applied anywhere in the nature), to conduct students in their scientific research, to advice them, to be involved in all academic activities and committee services, to enjoy teaching.

The first day of school can be more mathematically recreational. Ask the students: What do you like in mathematics, and what don't you like?

Tell them mathematics jokes, games, proofs with mistakes (to be found!), stories about mathematicians' lives, connections between mathematics and … opposite fields, such as: arts, music, literature, poetry, foreign languages, etc.

**2) What is the students' evaluation of you as instructor (negative opinions)?**

- do not be too nice in the classroom (because some students take advantage of that matter and waste their and the classmates time)
- be more strict and respond firmly
- don't say: "this is easy, you should know this" because one discourages students to ask questions
- attendance policy to be clear
- grammar skills, and listening skills
- patience with students
- allow students to help each other when they don't understand me
- clear English
- sometimes there is not enough time to cover all material
- to self-study the material and solve a lot of unassigned problems
- to talk louder to the class; to be more oriented towards the students and not to the board/self
- to understand student's questions
- to take off points if the homework problems are wrong, instead of just giving points for trying
- to challenge students in learning
- to give examples of harder problems on the board
- to enjoy teaching (smile, joke?)
- your methods should help students learning.

**3) What is the colleges' and university's mission and role in the society?**

- to ensure that all students served by the college learn the skills, knowledge, behaviors, and attitudes necessary for productive living in a changing, democratic, multicultural society.

**4) How do you see the future math teaching (new techniques)?**

- teaching on-line
- video-courses (with videotapes and tapes)
- teaching by using the internet
- teaching by using regular mail
- more electronic device tools in teaching (especially computers)
- interdisciplinary teaching
- self-teaching (helping students to teach themselves)
- more mathematics taught in connections with the social life (mathematical modeling)
- video conference style of teaching
- laboratory experiments

**5) What about "Creative Solutions"?**

- the focus of the program is on developing students' understanding of concepts and skills rather than "apparent understanding"
- students should be actively involved in problem-solving in new situations (creative solvers)
- non-routine problems should occur regularly in the student homework
- textbooks shall facilitate active involvement of students in the discovery of mathematical ideas
- students should make conjectures and guesses, experiment and formulate hypotheses and seek meaning
- the instructor should not let the teaching of mathematics degenerate into mechanical manipulation without thought
- to teach students how to think, how to investigate a problem, how to do research in their own, how to solve a problem for which no method or solution has been provided
- homework assignments should draw the student's attention to underlying concepts
- to do a cognitive guided instruction
- to solve non-routine problems, multi-step problems
- to use a step-by-step procedures for problem solving
- to integrate tradition with modern style teaching
- to emphasize the universality of mathematics
- to express mathematical ideas in a variety of ways
- to show students how to write mathematically, and how to read a mathematical text
- interpretations of solutions
- using MINITAB graphics to teach statistics (on the computer)
- tutorial programs on the computer
- developing manageable assessment procedures
- experimental teaching methods
- to motivate students to work and learn
- to stimulate mathematical reasoning

- to incorporate "real life" scenarios in the teacher training programs
- homo faber + homo sapiens are inseparable (Antonio Gramsci, Italian philosopher)
- to improve the critical thinking and reasoning skills of the students
- to teach students how to extend a concept
- to move from easy to medium and complex problems (gradually)
- mathematics is learned by doing, not by watching
- the students should be dedicate to the school
- to become familiar with symbols, rules, algorithms, key words and definitions
- to visualize mathematics notions
- to use computer generated patterns
- to use various problem solving strategies such as:
    - perseverance
    - achievement motivation
    - role model
    - confidence
    - flexible thinking
    - fresh ideas
    - different approaches
    - different data
- to use experimental teaching methods
- function plotters or computer algebra systems
- computer based learning
- software development
- grant proposal writing
- innovative pedagogy
- to use multi-representational strategies
- to try experimental tools
- to develop discussion groups
- symbol manipulation rules
- to solve template problems
- to do laboratory-based courses
- to think analytically
- to picture ourselves as teachers, or as students
- to use computer-generated patterns and new software tools
- to give to students educational and psychological tests to determine if any of them needs special education (for handicapped or gifted students) - American Association on mental Deficiency measures it.

6) **How to diminish the computer anxiety?**

In order to diminish the computer anxiety, a teacher needs to develop in students:
- positive attitudes towards appropriate computer usage
- feeling of confidence in use of computers
- feeling of comfort with computers

- acceptance of computers as a problem-solving tool
- willingness to use a computer for tasks
- attitude of responsibility for ethical use of the computer
- attitude that computers are not responsible for "errors"
- free of fear and intimidation of computers (the students' anxiety towards computer diminishes as their knowledge about computers increases)
- only after an algorithm is completely understood it is appropriate to rely on the computer to perform it
- computers help to remove the tedium of time – consuming calculations:
  - enable the students to consolidate the learning of the concepts and algorithms in math; the computer session is held at the end of the course when all the lectures and tutorials have been completed
  - to stimulate real world phenomena
  - all students should learn to use calculators
  - mathematics is easier if a calculator is used to solve problems
  - the calculator use is permissible on homework
  - using calculators makes students better problem solvers
  - the calculators make mathematics fun
  - using the calculator will make students try harder
  - the students should be able to
    a. assemble and start a computer
    b. understand the major parts of a computer
    c. use a variety of educational software
    d. distinguish the major instructional methodologies
    e. use word processor, database and spreadsheet programs
    f. attach and use a printer, peripherals, and lab probes
    g. use telecommunications networking
    h. use hypermedia technology
  - an instructor helps students to help themselves (it's interesting to study the epistemology of experience)

In the future the technology's role will increase due to the new kind of teaching: distance learning (internet, audio-visuals, etc.).

The technology is beneficial because the students do not waste time graphing functions anymore, but focusing on their interpretations.

7) **Describe your experience teaching developmental mathematics including course names, semester taught and methods and techniques used**.

In my teaching career of more than ten years experience I taught a variety of developmental mathematics courses such as:
- Introductory Mathematics: fall 1988, 1989; spring 1989, 1990. Methods: problem solving participation in the class, small groups work, guest speakers, discussions, student planning of assignments (to compose themselves problems of different styles and solve them by many methods), editing mathematics problem solutions (there are students who know how to solve a problem, but they are not able to write correctly and completely their proof mathematically), mathematics applied to real world problems (project), research work (how mathematics is used in a job), recreational mathematics approaches (logical games, jokes), etc.
- Pre-algebra: Spring 1982
- Algebra, Elem. Geometry: 1981, 1989.

**8) Briefly, describe your philosophy of teaching mathematics. Describe the application of this philosophy to a particular concept in a developmental mathematics course you have taught.**

- My teaching philosophy is "concept centered" as well as "problem solving directed". Makarenko said that everything can be taught to everybody if it's done at his/her level of knowledge". This focuses on promoting a student friendly environment where I not only lecture to provide to the student a knowledge base centering on concepts, but I also encourage peer mentoring with groups work to facilitate problem solving. It is my firm conviction that a student's perception, reasoning, and cognition can be strengthened with the application of both traditional and Alternative Learning Techniques and Student Interactive Activities.
- In my Introductory Mathematics course I taught about linear equations:
  - First I had to introduce the concept of variable, and then define the concept of equation; afterwards, tell to the students why the equation is called linear, how the linear equation is used in the real world, its importance in the every day's life;
  - Second, I gave students an example of solving a linear equation on the board, showing to them different methods, I classified them into consistent and inconsistent.

**9) Describe how you keep current with trends in mathematics instruction and give one example of how you have integrated such a trend into the classroom.**

- I keep current with trends in mathematics instruction reading journals such as" "Journal for Research in Mathematical Education", "Mathematics Teacher" (published by the National Council of Teachers of Mathematics, Reston, VA). "Journal of Computers in Mathematics and Science Teaching", "For the learning of mathematics", "Mathematics Teaching" (U. K.), "International Journal of Mathematical Education in Science and Technology", and participating with papers to the educational congresses, such as: The Fifth Conference on Teaching of Mathematics (Cambridge, June 21-22, 1996), etc. Example: Inter-subjectivity in Mathematics: teaching to everybody at his/her level of understanding.

10) **Describe your experience integrating technology into teaching mathematics. Provide specific examples of ways you have used technology in the mathematics classroom.**

- I use graphic calculators (TI-85) in teaching Intermediate Algebra; for example: programming it to solve a quadratic equation (in all 3 cases, when D is >, =, or < 0).
- I used various software packages of mathematics on IBM-PC or compatibles, such as: MPP, MAPLE, UA, etc. to give the students different approaches; for example in teaching Differential Equations I used MPP for solving a differential equation by Euler's method, changing many times the initial conditions, and graphing the solutions.

11) **Describe your knowledge and/or experience as related to your ability to prepare classroom materials.**

The classroom materials that I use: handouts, different color markers, geometric instruments, take-home projects, course notes, group projects, teaching outline, calculators, graphic calculators, PC, projectors, books, journals, etc.

12) **Describe the essential characteristics of an effective mathematics curriculum.**

- To develop courses and programs that support the College's vision of an educated person and a commitment to education as a lifelong process;
- To provide educational experiences designed to facilitate the individual's progress towards personal, academic, and work-based goals;
- To encourage the development of individual ideas and insights and acquisition of knowledge and skills that together result in an appreciation of cultural diversity and a quest for further discovery;
- To respond to the changing educational, social, and technological needs of current and prospective students and community employers.

13) **Provide specific examples of how you have and/or how you would develop and evaluate mathematics curriculum.**

In order to develop a mathematics curriculum:
I identify unmet student need, faculty interest in a new area, and requests from employers' recommendations of advisory committee, results of program review, university curriculum development.
Criteria for evaluation of a mathematics curriculum:
- course/program is educationally sound and positively affects course/program offerings within district; course does not necessarily duplicate existing course or course content in other disciplines offered throughout the district;

- development or modification of course/program does not adversely impact existing courses/programs offered throughout the district by competing for students and resources;
- course/program is compatible with the mission of the college.

**14) Describe your experience, education and training that has provided you with the knowledge and the ability to asses student achievement in mathematics.**

Courses I studied: History of Education, Introduction to Education, Philosophy, Child and Adolescent Psychology, Educational Psychology, General Psychology, Methods of Teaching Mathematics, Analysis of Teaching and Research, Instructional Design and Evaluation, Learning Skills Theory, Historical /Philosophical/ Social Education, Teaching Practice.

I taught mathematics in many countries, for many years, using various student assessments.

**15) Provide specific examples of ways you would asses student achievement in mathematics.**

I asses students by: tests in the Testing Center, quizzes in the classroom, homework, class participation (either solving problems on the board, or giving good answers for my questions), extra work (voluntarily), take home exams, research, projects, frequency.

Normally a test contains 10 problems, total being 100 points. For each homework I give 5 points, same for each extra work, for each class participation. For more than 3 absences I subtract points (one point for each absence), and later I withdraw the student.

Take home exams, quizzes, and research projects have the worth of a test.

Finally, I compute the average (my students know to assess themselves according to these rules, explained in the class and written in th syllabus).

**16) This question is about motivating a typical community college class of students, which is very diverse.**

- What kinds of students are you likely to have in such a class?
  Students of different races. genders, religions, ages, cultures, national origin, levels of preparedness, with or without physical or mental handicaps.
- How would you teach then?
  Catching their common interest, tutoring on a one-to one bases students after class (according to each individual level of preparedness, knowledge), working differentially with categories of students on groups, being a resource to all students, using multi representational strategies, motivating and making them dedicate to the study, finding common factors of the class. Varying teaching styles to respond to various student learning styles.

**17) Given the fact that the community college philosophy encourages faculty members to contribute to the campus, the college, and the community, provide examples of how you have and/or would contribute to the campus, the college, and the community.**

I have contributed to the college by:
- being an Associate Editor of the college (East Campus) "Math Power" journal;
- donating books, journals o the college (East Campus) Library;
- volunteering to help organizing the AMATYC math competitions (I have such experience from Romania and Morocco);
- representing the college at National/International Conferences on Mathematical and Educational Topics (for example at Bloomsburg University, PA, Nov. 14, 1995);
- publishing papers, and therefore making free publicity for the college;

I would contribute to the college by:
- organizing a math Club for interested students;
- cooperating with m fellow colleagues on educational projects sponsored by various foundations: National Science Foundation, Fulbright, Guggenheim;
- Socializing with m fellow colleagues to diverse activities needed to the college.
- Being a liaison between the College and University in order to frequently update the University math software and documentation (public property, reach done will a grant from NSF).

**18) Describe your experience within the last three years in teaching calculus for science and engineering majors and/or survey calculus at a post secondary level.**

I have taught Calculus I, II, III, in many countries. I have insisted on solving most creatively problems in calculus, because most of them are open-ended (they have more than one correct answer or approach); sometimes, solving a problem relies on common sense ideas that are not stated in the problem. The fundamental basis of the Calculus class is what graphs symbolize, not how to draw them.
Using calculators or computers the students got reasonable approximation of a solution, which was usually just useful as an exact one.

**19) Reform calculus a significant issue in math education today. Describe your thoughts on the strengths and weaknesses of reform versus traditional calculus and indicate which form of calculus you would refer to teach.**

Of course, I prefer to teach the Harvard Calculus, because it gives the students the skills t read graphs and think graphically, to read tables and think numerically, and to apply these skills along with their algebraic skills to modeling the real world (The Rule of Three); and Harvard Calculus also states that formal mathematical theory evolves from investigations of practical problems (The Way of Archimedes).

Weaknesses: the students might rely too much on calculators or computers ("the machines will think for us!"), forgetting to graph, solve, computer.

**20) Describe your experience in curriculum development including course development, textbook or lab manual development, and development of alternative or innovative instructional methods.**

I have developed a course of Calculus I, wrote and published a textbook of Calculus I for students, associated with various problems and solutions on the topic.

Concerning the alternative instructional methods, I'm studying and developing The Inter-subjectivity Method of Teaching in Mathematics (inspired by some articles from "Journal for Research in Mathematical Education" and "International Journal of Mathematical Education in Science and Technology").

**21) Describe your education and/or experiences that would demonstrate your ability to proactively interact with and effectively teach students from each of the following: different races, cultures, ages, genders, and levels of preparedness. Provide examples of your interaction with and teaching of students from each of these groups.**

I have taught mathematics in many countries: Romania (Europe), Morocco (Africa), Turkey (Asia), and U.S.A. Therefore, I am accustomed to work with a diverse student population. More, each country had its educational rules, methods, styles, curriculum missions – including courses, programs, textbooks, mathematics student competitions, etc. that I have acquired a very large experience. I like to work in a multi-cultural environment teaching in many languages, styles (according to the students' characteristics), being in touch with various professors around the world, knowing many cultural habits.

**22) Describe your professional development activities that help you stay in the field of mathematics. Give your best example of how you have integrated one thing into the classroom that came out of your professional development activities.**

I subscribe to math journals, such as: "College Mathematics Journal", as a member of the Mathematical Association of America, and often go to the University Libraries, Science Section, to consult various publications.

I keep in touch with mathematicians and educators from all over, exchanging math papers and ideas, or meeting them at Conferences or Congresses of math or education.

Studying about "intersurjectivity" in teaching, I got the idea of working differentially with my students, distributing them in groups of low level, medium level, high level according to their knowledge, and therefore assigning them appropriately special projects.

**23)**
   a) **What are the most important personal and academic characteristics of a teacher?**
   b) **At the end of your first year of district employment how will you determine whether or not you have bee successful?**
   c) **What are the greatest challenges in public education today?**
   d) **What do you want your students to learn?**

- a) To be very good professional in his/hers field, improving his/her skills permanently. To be dedicated to his/her work. To love the students and understand their psychology. To be a very good educator. To prepare every lesson (its objectives). To provide attractive and interesting lessons.
- b) Regarding the level of the class (the knowledge in math), the student's grades, even their hobby for mathematics (or at least their interest).
- c) To give the students a necessary luggage of knowledge and enough education such that they are able to fend themselves in our society (they are prepared very well for the future)
- d) To think. Brainstorm. To solve not only mathematical problems, but also life problems.

**24) What do you want to accomplish as a teacher?**

I like to get well prepared students with good behaviors.

**25) How will you go about finding our students' attitudes and feelings about your class?**

I'll try to talk with every student to find out their opinions, difficulties, and attitudes towards the teacher. Then, I'll try to adapt myself to the class level of knowledge and be agreeable to the students. Besides that, I'll try to approach them in extracurricular activities soccer, tennis, chess, creative art and literature using mathematical algorithms/methods, improving my Spanish language.

**26) An experienced teacher offers you following advice: "When you are teaching, be sure to command the respect of your students immediately and all will go well". How do you feel about this?**

I agree that in a good lesson the students should respect their teacher, and reciprocally. But the respect should not be "commanded", but earned. The teacher should not hurt the students by his/her words.

**27) How do you go about deciding what it is that should be taught in your class?**

I follow the school's plan, the mathematics text book; the school's governing board directions. I talk with other mathematics teachers asking their opinions.

**28) A parent comes to you and complains that what you are teaching his child is irrelevant to the child needs. How will you respond?**

I try to find out what the parent wants, what his needs are. Then, maybe I have to change my teaching style. I respond that irrelevant subjects of today will be relevant subjects of tomorrow.

**29) What do you think will (does) provide you the greatest pleasure in teaching?**

When students understand what I'm teaching, and they know how to use what I taught them in real life.

**30) When you have some free time, what do you enjoy doing the most?**

I try to improve my mathematical skills (subscriptions to mathematical and educational journals). Teaching mathematics became a hobby for me!

**31) How do you go about finding what students are good at?**

I try to approach mathematics with what students are good at. For example: I tell them that mathematics are applied anywhere in the nature and society, therefore in arts, in music, in literature, etc. Therefore, we can find a tangential joint between two apparent distinct (opposite) interests.

**32) Would you rather try a lot of way-out teaching strategies or would you rather try to perfect the approaches, which work best for you? Explain your position.**

Both: the way-out teaching strategies combined with approaches to students. In each case the teacher should use the method/strategy that works better.

**33) Do you like to teach with an overall plan in mind for the year, or would you rather just teach some interesting things and let the process determine the results? Explain your position.**

Normally I like to teach with an overall plan in mind, but some times – according to the class level and feelings – I may use the second strategy.

**34) A student is doing poorly in your class. You talk to him/her, and he/she tells you that he/she considers you to be the poorest teacher he/she ever met. What would you do?**

I try to find out the opinions of other students about my teaching and to get a general opinion of the entire class. I give students a test with questions about my character, skills, style, teaching methods, etc. in order to find out my negative features and to correct/improve them by working hard.

**35) If there were absolutely no restrictions placed upon you, what would you most want to do in life?**

I would like to set up a school (of mathematics especially) for gifted and talented students with a mathematics club for preparing students for school competitions.

**36) How do you test what you teach?**

By written test, final exams, homework, class participation, special projects, extra homework, quizzes, and take-home exams.

**37) Do you have and follow a course outline? When would a variation from the outline be appropriate?**

- Yes, I follow a course outline.
- When I find out the students have gaps in their knowledge and, therefore they are not able to understand the next topic to be taught. Or new topics are needed (due to scientific research or related other disciplines).

**38) Is student attendance important for your course? Why or why not? What are the student responsibilities necessary for success in your class?**

- Yes.
- If they miss many courses they will have difficulties to understand the others, because mathematics is like a chain.
- To work in the classroom, to pay attention and ask questions, to do independent study at home too.

**39) Describe your turnaround time for returning graded test and assignments.**

I normally grade the tests over the weekends. The same I do for all other assignments.

**40) Are you satisfied with the present textbooks? Why or why not?**

- Yes.
- Because they gives the students the main ideas necessary in the technical world.

**41) Describe some of the supplemental materials you might use for this course.**

- Personal computers with DERIVE software package.
- T!-92 and an overhead projector.
- Tables of Laplace Transforms.
- Various handouts.

**42) Describe your method of student record keeping.**

- I keep track of: absences, homework, tests' grades, final exam's grade, class participation.

**43) Describe how you assist or refer students who need remediation.**

- I advice them to go to College Tutoring Center.
- I encourage them to ask questions in the classroom, to work in groups with the better students, to contact me before or after class.

**44) What is your procedure for giving students feedback on their learning progress?**

- By the test grades
- By the work they are doing in the classroom.

**45) How do you monitor your evaluation methods so that they are both fair and constructive?**

- My students are motivated to work and improve their grades by doing extra (home) work.
- I compare my evaluation methods with other instructors'.
- I also fell when a student masters or not a subject.

**46) Describe your relationship with your colleagues.**

- I share information, journals, books, samples of tests, etc. with them
- Good communication.

**47) What procedures do you use to motivate students?**

- Give them a chance to improve their grades.
- Telling them that if they don't learn a subject in mathematics, they would not understand the others (because mathematics is cyclic and linear).

**48) Are you acquainted with district and campus policies and procedures? Do you have any problems with any of the policies and procedures?**

- I always try to adjust myself to each campus' policy

**49) What mathematical education topic are you working in?**

- I am studying the radical constructivism (Jean Piaget) and social constructivism (Vygotsky: to place communication and social life at the

center of meaning – making), the inter-subjectivity in mathematics, the meta-knowledge, the assessment standards. Learning and teaching are processes of acculturation.